\documentclass[12pt]{amsart}

\usepackage{amsmath}
\usepackage{latexsym}
\usepackage{amssymb}
\usepackage{amsfonts}
\usepackage{graphics}
\usepackage{color}

\renewcommand{\proof}{\par\noindent{\it Proof.\ \ }}
\def\qed{\ifmmode\square\else\nolinebreak\hfill
$\Box$\fi\par\vskip12pt}

\def\div{\,\big|\,} \def\notdiv{{\,\not\big|\,}}

    \def\ZZ{\mathbb Z}

\def\mod{{\sf mod~}}

\def\soc{{\sf soc}} \def\rad{{\sf rad}}

\def\N{{\bf N}}

\def\Ga{{\it \Gamma}}

\def\Cay{{\sf Cay}} \def\Cos{{\sf Cos}} 
\def\Aut{{\sf Aut}} \def\Inn{{\sf Inn}} \def\Out{{\sf Out}}

\def\D{{\sf D}} 
\def\S{{\sf S}} 
 \def\M{{\sf M}}
\def\Mult{{\sf Mult}}

\def\a{\alpha} \def\b{\beta} \def\g{\gamma}

\def\POmega{{\sf P\Omega}}
\def\Sp{{\sf Sp}}\def\PSp{{\sf PSp}}

\def\AGammaL{{\sf A\Gamma L}}
 
\def\A{{\sf A}} 

\def\PSL{{\sf PSL}}\def\PGL{{\sf PGL}}
\def\GL{{\sf GL}}

\def\AGL{{\sf AGL}}

\def\AGammaL{{\sf A\Gamma L}}

 \def\PSU{{\sf PSU}} \def\SU{{\sf SU}} 
  \def\D{{\sf D}}

\def\mz{{\mathbb Z}}

\newtheorem{theorem}{Theorem}[section]%
\newtheorem{lemma}[theorem]{Lemma}%
\newtheorem{proposition}[theorem]{Proposition}%
\newtheorem{problem}[theorem]{Problem}%

\begin{document}

\title[ Arc-transitive Cayley graphs on nonabelian simple groups]
{Arc-transitive Cayley graphs on nonabelian simple groups with prime valency}%
\author{Fu-Gang Yin, Yan-Quan Feng, Jin-Xin Zhou, Shan-Shan Chen}%
\address{Department of Mathematics, Beijing Jiaotong University, Beijing, 100044, China}
 \email{181181010@bjtu.edu.cn (F.-G. Yin), yqfeng@bjtu.edu.cn (Y.-Q. Feng), \linebreak
  jxzhou@bjtu.edu.cn (J.-X. Zhou), 18121630@bjtu.edu.cn (S.-S. Cheng) } %

 \thanks{This work was supported by the National Natural Science Foundation of China (11731002,11671030) and by the 111 Project of China (B16002).}
 \subjclass{05C25, 20B25}%


\begin{abstract}
In 2011, Fang et al. in (J. Combin. Theory A 118 (2011) 1039-1051) posed the following problem: {\em Classify non-normal locally primitive Cayley graphs of finite simple groups of valency $d$, where either $d\leq 20$ or $d$ is a prime number}. The only case for which the complete solution of this problem is known is of $d=3$. Except this, a lot of efforts have been made to attack this problem by considering the following problem: {\em Characterize finite nonabelian simple groups which admit non-normal locally primitive Cayley graphs of certain valency $d\geq4$.} Even for this problem, it was only solved for the cases when either $d\leq 5$ or $d=7$ and the vertex stabilizer is solvable. In this paper, we make crucial progress towards the above problems by completely solving the second problem for the case when $d\geq 11$ is a prime and the vertex stabilizer is solvable.

\end{abstract}
\maketitle
 \qquad {\textsc k}{\scriptsize \textsc {eywords.}} {\footnotesize  Cayley graph, simple group, arc-transitive graph.}

%
\maketitle

\section{Introduction}
Throughout this paper, graphs are assumed to be finite undirected graphs without loops and multiple edges, and groups are assumed to be finite. Let $G$ be a permutation group on a set $\Omega$, and let $\a\in \Omega$. Denote by $G_{\a}$ {\it the stabiliser} of $\a$ in $G$, that is, the subgroup of $G$ fixing the point $\a$. The group $G$ is {\it semiregular} if $G_{\a}=1$ for every $\a\in\Omega$, and {\it regular} if $G$ is transitive and semiregular.

For a graph $\Ga$, denote by $V(\Ga)$, $E(\Ga)$ and $\Aut(\Ga)$ its vertex set, edge set and full automorphism group, respectively. For a vertex $v \in V(\Ga)$, let $\Ga(v)$ be the neighbourhood of $v$ in $\Ga$. An {\it $s$-arc} in $\Ga$ is an ordered $(s+1)$-tuple $(v_0,v_1,...,v_s)$ of vertices of $\Ga$ such that $v_{i-1}$ is adjacent to $v_i$ for $1\leq i\leq s$, and $v_{i-1}\neq v_{i+1}$ for $1\leq i<s$. A graph $\Ga$, with $G\leq \Aut(\Ga)$, is said to be {\it $(G,s)$-arc-transitive} or  {\it $G$-regular} if $G$ is transitive on the $s$-arc set of $\Ga$ or $G$ is regular on the vertex set $V(\Ga)$ of $\Ga$, respectively. For short, a $1$-arc means an {\it arc}, and  $(G,1)$-arc-transitive means {\it $G$-arc-transitive}. If a graph $\Ga$ is $G$-regular, then $\Ga$ is also called a {\em Cayley graph} of $G$, and the Cayley graph is \emph{normal} if $G$ is normal in $\Aut(\Ga)$.
A graph $\Ga$ is said to be {\it $s$-arc-transitive} if it is $(\Aut(\Ga), s)$-arc-transitive. In particular, $0$-arc-transitive is {\it vertex-transitive}, and $1$-arc-transitive is {\it arc-transitive} or {\it symmetric}.

A fair amount of work have been done on symmetric Cayley graphs on non-abeian simple groups in the literature. One of the remarkable achievements in this research field is the complete classification of cubic non-normal symmetric Cayley graphs of non-abelian simple groups, and it turns out that up to isomorphism, there are only two cubic non-normal symmetric Cayley graphs of non-abelian simple groups which are both cubic $5$-arc-transitive Cayley graphs on $\A_{47}$ (see \cite{Li,XFWX1,XFWX2}). Recall that a graph $\Ga$ is called {\em locally primitive} if for any $v\in V(\Ga)$, the stabilizer $\Aut(\Ga)_v$ of $v$ in $\Aut(\Ga)$ is primitive on $\Ga(v)$. In view of the fact that every cubic symmetric graph is locally primitive, a natural question arises: What can we say about locally primitive non-normal symmetric Cayley graphs of non-abelian simple groups?

On locally primitive graphs, Weiss~\cite{Weiss} conjectured that there
is a function $f$ defined on the positive integers such that, whenever $\Ga$ is a $G$-vertex-transitive locally primitive graph of valency $d$ with $G\leq \Aut(\Ga)$ then, for any vertex $v\in V(\Ga)$, $|G_v|\leq f (d)$. By Conder et al.~\cite{Conder-Li-Preger-2000}, Weiss conjecture is true for vertex-transitive locally primitive $d$-valent graphs if $d\leq 20$ or $d$ is a prime number, and by Spiga~\cite{Spiga}, Weiss conjecture is also true if the restriction $G^{\Ga(v)}$ of $G$ on $\Ga(v)$ contains an abelian regular subgroup, that is, of affine type. In 2007, Fang et al.~\cite[Theorem~1.1]{FMW} shown that for any valency $d$ for which the Weiss conjecture holds, all but finitely many locally primitive Cayley graphs of valency $d$ on the finite nonabelian simple groups are normal, and based on this, the following problem was proposed:

\begin{problem}{\rm\cite[Problem~1.2]{FMW}}\label{prob1}
Classify non-normal locally primitive Cayley graphs of finite simple groups of valency $d$, where either $d\leq 20$ or $d$ is a prime number.
\end{problem}

As mentioned above, this problem has been completely solved by Li {\em et al.} for the case when $d=3$. For the case when $d\geq 4$, however, it is quite difficult to give a complete solution of Problem~\ref{prob1}. Because of this, researchers have focused on the following slightly easier problem.

\begin{problem}\label{prob2}
Characterize finite nonabelian simple groups which admit non-normal locally primitive Cayley graphs of certain valency $d\geq 4$.
\end{problem}

Clearly, a tetravalent graph is locally primitive if and only if the graph is $2$-arc-transitive. In 2004, Fang {\em et al}~\cite{FLX} proved that except 22 groups given in \cite[Table~1]{FLX}, every tetravalent $2$-arc-transitive Cayley graph $\Ga$ of a non-abelian simple group $G$ is normal, and based on this, in 2018, Du and Feng~\cite{DF} proved that there are exactly $7$ non-abelian simple groups which admit at least one non-normal $2$-arc-transitive Cayley graph, thus giving a complete solution of Problem~\ref{prob2} for the case when $d=4$.

There are also some partial solutions of Problem~\ref{prob2} for the case when $d$ is a prime number. It is easy to see that a graph with prime valency is locally primitive if and only if it is symmetric. Fang {\em et al}~in \cite{FMW} constructed an infinite family of $p$-valent non-normal symmetric Cayley graphs of the alternating groups for all prime $p\geq 5$, and using a result in \cite{FPW} on the automorphism groups of Cayley graphs of non-abelian simple groups, they also gave all possible candidates of finite nonabelian simple groups which might have a pentavalent non-normal symmetric Cayley graph. This was recently improved by Du {\em et al}~\cite{DFZ} by proving that there are only $13$ finite nonabelian simple groups which admit a pentavalent non-normal symmetric Cayley graph.

More recently, Pan {\em et al}~\cite{PYL} considered Problem~\ref{prob2} for the case when $d=7$, and they proved that for a $7$-valent Cayley graph $\Ga$ of a non-abelian simple group $G$ with solvable vertex stabilizer, either $\Ga$ is normal, or $\Aut(\Ga)$ has a normal arc-transitive non-abelian simple subgroup $T$ such that $G<T$ and $(G,T)=(\A_6,\A_7)$, $(\A_{20},\A_{21})$, $(\A_{62},\A_{63})$ or $(\A_{83},\A_{84})$, and for each of these $4$ pairs $(G,T)$, there do exist a $7$-valent $G$-regular $T$-arc-transitive graph.

In this paper, we shall prove the following theorem which generalizes the result in \cite{PYL} to all prime valent cases, and hence gives a solution of Problem~\ref{prob2} for the case when $d$ is a prime and the vertex-stabilizer is solvable.

\begin{theorem}\label{th1}
Let $G$ be a non-abelian simple group and $\Ga$ a connected arc-transitive Cayley graph of $G$ with prime valency $p \geq 11$. If $\Aut(\Ga)_v$ is solvable for $v \in V(\Ga)$, then either $G \unlhd \Aut(\Ga)$, or $\Aut(\Ga)$ has a normal subgroup $T$ with $G<T$ such that $\Ga$ is $T$-arc-transitive and $(G, T, p)$ is one of the following four triples:
\[(\A_5, \PSL(2,11), 11), (\A_5, \PSL(2,29), 29),  (\M_{22}, \M_{23}, 23), (\A_{n-1}, \A_n, p),\]
where $n=pk\ell$ with $k\div \ell$ and $\ell\div (p-1)$, and $k$ and $\ell$ have the same parity.
\end{theorem}

Conversely, we show that all the first three triples as well as the fourth triple in case of $n=p$ can happen.

\begin{theorem}\label{th2}
Use the same notation as {\rm Theorem~\ref{th1}}. If $(G, T, p)$ is one of the following triples:
\[(\A_5, \PSL(2,11), 11), (\A_5, \PSL(2,29), 29),  (\M_{22}, \M_{23}, 23), (\A_{p-1}, \A_p, p),\]
then there exists a $p$-valent symmetric Cayley graph $\Ga$ of $G$ such that $\Aut(\Ga)_v$ is solvable for some $v\in V(\Ga)$.
\end{theorem}

Let $p$ be a prime and $\ell,k$ integers with $k\div \ell$ and $\ell\div (p-1)$ such that $k$ and $\ell$ have the same parity.
The triple $(p,\ell,k)$ is called {\em conceivable} if there exists an arc-transitive Cayley graph of the alternating group $\A_{pk\ell-1}$ with valency $p$ and its automorphism group has solvable vertex stabilizer. We have been unable to determine all the conceivable triples $(p,\ell,k)$, and we would like to leave it as an open problem for future research.

\begin{problem}
Determine conceivable triples $(p,\ell,k)$.
\end{problem}

By Theorem~\ref{th2}, $(p,1,1)$ is conceivable for each prime $p\geq 5$, and by \cite{DFZ}, $(5,4,2)$ is conceivable, but not $(5,2,2)$. For the case $p=7$, it was shown in \cite{PYL} that $(7,1,1)$, $(7,3,1)$, $(7,3,3)$ and $(7,6,2)$ are the only conceivable triples.

The paper is organized as follows. In Section 2 we introduce some preliminary results on nonabelian simple groups and arc-transitive graphs with prime valency. Then we prove Theorem \ref{th1} in Section~\ref{proofth1} and Theorem~\ref{th2} in Section~\ref{nonnormalexample}.


\section{Preliminary}

In this section, we introduce some preliminary results that will be used latter.

For a positive integer $n$, we use $\ZZ_{n}$ to denote the cyclic group of order $n$.
For a group $G$ and a subgroup $H$ of $G$, denote by $N_G (H)$ and $C_G (H)$ the normalizer and the centralizer of $H$ in $G$ respectively.
Given two groups $N$ and $H$, denote by $N \times H$ the direct product of $N$ and $H$, by $N.H$ an extension of $N$ by $ H$, and if such an extension is split, then we write $N : H$ instead of $N . H$.

 The following proposition is an exercise in Dixon and Mortimer's textbook \cite[p.49]{Dixon}.

\begin{proposition} \label{lmsylown}
Let  $n$  be  a positive integer  and  $p$  a prime. Let $p^{\nu(n)} $ be the largest  power of  $p$  which divides  $n!$. Then $\nu(n)=\sum_{i=1}\lfloor \frac{n}{p^{i}} \rfloor< \frac{n}{p-1}$.
\end{proposition}

The next proposition is called the {\it Frattini argument } on transitive permutation group, and we refer to \cite[p.9]{Dixon}.

\begin{proposition}\label{Frattiniaugrment}
 Let $G$ be a transitive permutation group on $\Omega$, $H$ a subgroup of $G$ and $v \in \Omega$. Then $H$ is transitive if and only if $G = HG_v $.
\end{proposition}

We denote by $\Aut(G)$ the automorphism group of a group $G$, and by $\Inn(G)$ the inner automorphism group of $G$ consisting of these automorphisms of $G$ induced by all element of $G$ by conjugation on $G$. Then $\Inn(G)$ is normal in $\Aut(G)$, and the quotient group $\Aut(G)/\Inn(G)$ is called the {\em outer automorphism} of $G$, denoted by $\Out(G)$. The following proposition is a direct consequence of the classification of finite simple groups (see \cite[Table 5.1.A-C]{K-Lie} for example).

\begin{proposition} \label{thoutT}
Let $T$ be a finite non-abelian simple group. Then $\Out(T)$ is solvable.
\end{proposition}

Let $G$ and $E$ be two groups. We call an extension $E$ of $G$ by $N$ a {\em central extension} of $G$ if $E$ has a central subgroup $N$ such that $E/N\cong G$, and if further $E$ is perfect, that is, the derived group $E'$ equals to $E$, we call $E$ a {\em covering group} of $G$. A covering group $E$ of $G$ is called a {\em double cover} if $|E|=2|G|$. Schur~\cite{Schur} proved that for every non-abelian simple group $G$ there is a unique maximal covering group $M$ such that every covering group of $G$ is a factor group of $M$ (see \cite[Kapitel V, \S23]{Huppert}). This group $M$ is called the {\em full covering group} of $G$, and the center of $M$ is the {\em Schur multiplier} of $G$, denoted by $\Mult(G)$. By Kleidman and Liebeck~\cite[Theorem 5.1.4]{K-Lie} and Du {\em et al}~\cite[Proposition 2.6]{DFZ}, we have the following proposition.

\begin{proposition}\label{lmcoveran} $\Mult(\A_{n})=\ZZ_2$ with $n \geq 8$. For $n\geq 5$, $\A_n$ has a unique double cover $2.\A_n$, and for $n\geq 7$, all subgroups of index $n$ of $2.\A_n$ are isomorphic to $2.\A_{n-1}$.
\end{proposition}

By Kleidman and Liebeck~\cite[Proposition 5.3.7]{K-Lie}, we have the following proposition.

\begin{proposition} \label{lm-An}
Let $r$ be a prime power and $f$ a positive integer.  If $\A_n \leq \GL(f,r) $ with $n \geq 9$, then $f \geq n-2$.
\end{proposition}

Let $\Ga$ be a connected graph and $G$ a group of automorphisms of $\Ga$. For $v\in V(\Ga)$, denote by $G_v^{\Ga(v)}$ the induced permutation group of the natural action of $G_v$ on the neighbourhood $\Ga(v)$. Let $G_v^*$ be the subgroup of $G_v$ fixing every vertex in $\Ga(v)$. Then $G_v^*$ is the kernel of the natural action of $G_v$ on $\Ga(v)$, and hence $G_v/G_v^*\cong G_v^{\Ga(v)}$. By the connectivity of $\Ga$, there exists a path $v=v_0,v_1,v_2,\cdots,v_m$ such that $G_{v_0v_1\cdots v_m}^*:=G_{v_0}^*\cap G_{v_1}^*\cap\cdots \cap G_{v_m}^*=1$. Clearly,
$$1=G_{v_0v_1\cdots v_m}^*\unlhd G_{v_0v_1\cdots v_{m-1}}^*\unlhd \cdots \unlhd G_{v_0v_1}^*\unlhd G_{v_0}^*=G_v^*\unlhd G_v,$$
and for $0\leq i< m$, we have $G_{v_0v_1\cdots v_i}^*/G_{v_0v_1\cdots v_{i+1}}^*\cong (G_{v_0v_1\cdots v_i}^*)^{\Ga(v_{i+1})}$. Then we can easily obtain the following proposition, and this was known from a series of lectures given by Cai Heng Li in Peking University in 2013.

\begin{proposition}\label{localgloblestabilizers}
Let $\Ga$ be a connected graph and let $G$ be a vertex-transitive group of automorphisms of $\Ga$. Then $G_v$ is nonsolvable if and only if $G_v^{\Ga(v)}$ is nonsolvable.
\end{proposition}

For self-containing, we give a short proof of the following proposition, which is mainly owed to an anonymous referee (also see \cite{Guo} for another proof).

\begin{proposition} \label{solvable-stabilizer}
Let $\Ga$ be a connected $G$-arc-transitive graph of prime valency $p\geq 5$, and let $(u,v)$ be an arc of $\Ga$. Assume that $G_v$ is solvable. Then $G_{uv}^*=1$ and $G_v\cong \mz_k \times (\mz_p:\mz_\ell)$ with $k\div \ell\div (p-1)$, where $\mz_p:\mz_\ell \leq \AGL(1,p)$.
\end{proposition}

\proof It follows from \cite{Weiss} that $G_{uv}^*= 1$. Let $P$ be a Sylow $p$-subgroup of $G_v$. Note that $G_v^{\Ga(v)}$ is a transitive solvable group of prime degree. By the Burnside Theorem (also see \cite[Theorem  3.5B]{Dixon}), $G_v/G_v^*\cong G_v^{\Ga(v)}\cong \mz_p:\mz_\ell \leq \AGL(1,p)$ with $\ell\div (p-1)$ and $G_{uv}/G_v^*\cong \mz_\ell$. In particular, $PG_v^*/G_v^*\unlhd G_v/G_v^*$, and so $PG_v^*\unlhd G_v$.  Since $G_u^*=G_u^*/G_{uv}^*=G_u^*/(G_u^*\cap G_v^*)\cong G_u^*G_v^*/G_v^*\leq G_{uv}/G_v^*\cong \mz_\ell$, we have $G_v^*\cong \mz_k$ with $k\div \ell$, and then $|G_v|=pk\ell$ with $k\div \ell\div (p-1)$. Since $G_{uv}=G_{uv}/G_{uv}^*= G_{uv}/(G_u^*\cap G_v^*)\lessapprox  G_{uv}/G_u^*\times  G_{uv}/ G_v^*\cong \mz_\ell\times \mz_\ell$, $G_{uv}$ is abelian of exponent $\ell$. Let $G_{uv}/G_v^*=\langle aG_v^*\rangle\cong \mz_\ell$. Then $\langle a\rangle\cong\mz_\ell$ and $\langle a\rangle\cap G_v^*=1$. It follows that $G_{uv}=\langle a\rangle \times G_v^*$.

Since $|G_v^*|=\ell\div (p-1)$, $PG_v^*$ has a unique Sylow $p$-subgroup $P$ and hence $PG_v^*=P\times G_v^*$. Then $P$ is characteristic in $PG_v^*$, and since $PG_v^*\unlhd G_v$, we have $P\unlhd G_v$. It follows that $G_v=P: G_{uv}=P: (\langle a\rangle\times G_v^*)=G_v^* \times (P: \langle a\rangle) \cong \mz_k \times (\mz_p: \mz_\ell)$. \qed

Taking normal quotient graphs is a useful method for studying arc-transitive graphs. Let $\Ga$ be an $X$-vertex-transitive graph, where $X \leq \Aut(\Ga)$ has an intransitive normal subgroup $N$. The {\it normal quotient graph } $\Ga_ N$ of $\Ga$ induced by $N$ is defined to be a graph with vertex set $\{ \alpha^N\ |\ \alpha \in V(\Ga) \}$, the set of all $N$-orbits on $V(\Ga)$, such that two vertices $B , C \in \{ \alpha^N\ |\ \alpha \in V(\Ga) \}$ are
adjacent if and only if some vertex in $B$ is adjacent in $\Ga$ to some vertex in $C$. If $\Ga$ and $\Ga_N$ have the same valency, then $\Ga$
is called a {\em normal cover } of  $\Ga_N$.
The following proposition is a special case of \cite[Lemma 2.5]{LP}, which slightly improves a remarkable result of Praeger \cite[Theorem 4.1]{Praeger}.

\begin{proposition}\label{lm-arc}
Let $\Ga $ be a connected $X$-arc-transitive graph of prime valency, with $X\leq \Aut(\Ga)$, and let $N ~ \unlhd~ X$ have at least three orbits on $V(\Ga)$. Then the following statements hold.
\begin{itemize}
\item[\rm (1)] $N$ is semi-regular on $V(\Ga)$, $X/N \leq \Aut(\Ga _N)$, $\Ga_N$ is a connected $X/N$-arc-transitive graph, and $\Ga$ is a normal cover of $\Ga_N$.

\item[\rm (2)] $X_v \cong (X/N)_{\Delta}$ for any $v \in V(\Ga)$ and $\Delta \in V(\Ga_N)$.
 \end{itemize}

\end{proposition}

\section{Proof of Theorem \ref{th1} \label{proofth1}}

Throughout this section we make the following assumption.

\medskip

\noindent{\bf Assumption:} $\Ga$ is a symmetric graph of prime valency $p \geq 11$ with $v\in V(\Ga)$, $\Aut(\Ga)_v$ is solvable, and $G\leq \Aut(\Ga)$ is a non-abelian simple group and transitive on $V(\Ga)$.

\medskip
The proof of the following lemma is straightforward, but will be used frequently latter.

\begin{lemma}\label{lm-homo}
Let $X=H:K$ be a transitive permutation group on $\Omega$. Let $w \in \Omega$. If $H $ is transitive, then $K$ is isomorphic to $X_w/H_w$.
\end{lemma}

\proof Since $H$ is transitive, $X=HX_w$ by Proposition~\ref{Frattiniaugrment}. So $K \cong X/H=HX_w/H
\cong X_w/(X_w \cap H)=X_w/H_w$. \qed

The product of all minimal normal subgroups of a group $X$ is called the {\em socle} of $X$, denoted by $\soc(X)$, and the largest normal solvable subgroup of $X$ is called the {\em radical} of $X$, denoted by $\rad(X)$.

\begin{lemma}\label{lm-R=1-trans} Let $G$, $\Ga$, $p$ and $v$ be as given in Assumption. Let $\Ga$ be $X$-arc-transitive with $G \leq X\leq \Aut(\Ga)$, and let $\rad(X)=1$. Then either $\soc(X)=G$, or $\Ga$ is $\soc(X)$-arc-transitive with $G<\soc(X)$ and one of the following holds:
\begin{itemize}
\item[\rm (1)] \ $(G,\soc(X))=(\A_{n-1}, \A_n)$ with $n\geq 6$, and $(\soc(X))_v$ is transitive on $\{ 1,2, \cdots,n \}$.
\item[\rm (2)]  $(G,\soc(X))=(\M_{22},\M_{23})$, and $(\soc(X))_v=\ZZ_{23}$.
\item[\rm (3)] \ $(G,\soc(X))=(\A_5,\PSL(2,11))$, and $(\soc(X))_v=\ZZ_{11}$.
\item[\rm (4)] \ $(G,\soc(X))=(\A_5,\PSL(2,29))$, and $(\soc(X))_v=\ZZ_{29}:\ZZ_{7}$.

\end{itemize}
In particular, $\Ga$ is a Cayley graph of $G$ for cases (2)-(4).

\end{lemma}
\proof Let $N$ be a minimal normal subgroup of $X$. Since $\rad(X)=1$, we have $N=T_1 \times \cdots \times T_d\cong T^d$ for a non-abelian simple group $T$. Write $K=NG$.

Assume that $G\unlhd X$. If $N\cap G=1$, applying Lemma \ref{lm-homo} with $K=G :N$ we have that $N \cong (K)_v/G_v$ is solvable, a contradiction. Therefore, $N\cap G\not=1$, forcing $G\leq N$, and since $G$ is normal, the minimality of $N$ implies $N = G$. By the arbitrariness of $N$, we have $\soc(X)=G$.

In what follows we assume that $G\ntrianglelefteq X$. If $\Ga$ is bipartite, then the transitivity of $G$ on $V(\Ga)$ implies that $G$ has a normal subgroup of index $2$, contradicting the simplicity of $G$. Thus, $\Ga$ is not bipartite. Therefore $N$ has either one or at least three orbits on $V(\Ga)$. We claim that the latter cannot occur.

We argue by contradiction and we suppose that $N$ has at least three orbits on $V(\Ga)$.
By Proposition \ref{lm-arc}, $N$ is semiregular on $V(\Ga)$, and so $|N|=|T|^d$ is a divisor of $|V(\Ga)|$. In particular, $|N| \mid |G|$. Since $N$ has at least three orbits, $|G|\geq 3|N|$ and hence $N \cap G=1$.

Consider the conjugate action of $G$ on $N$, and since $G$ is simple, the action is trivial or faithful. If it is trivial then $K=N\times G$, and by Lemma~\ref{lm-homo}, $N\cong K_v/G_v$ is solvable, a contradiction. It follows that the conjugate action of $G$ on $N$ is faithful, and hence we may assume $G\leq \Aut(N)$.

Note that $\Aut(N) \cong \Aut(T)^d:\S_d$. Set $M=\Aut(T)^d$ and $M_1=\Inn(N) \cong T^d$. Then $|M_1|=|N|$, $M_1\unlhd M$, $M\unlhd \Aut(N)$ and $M_1\unlhd \Aut(N)$. Clearly, $G\cap M_1=1$ as $|G|\geq 3|N|=3|M_1|$. If $G\cap M\not=1$ then $G\leq M$ and hence $G\cong G/(G\cap M_1)\cong GM_1/M_1\leq M/M_1\cong \Out(T)^d$, which is impossible because $\Out(T)$ is solvable by Propostion~\ref{thoutT}. This means that $G\cap M=1$, and therefore, $G\cong G/(G\cap M)\cong GM/M\leq \Aut(N)/M\cong \S_d$.
Recall that $|N|=|T|^d$ and $|N|\mid |G|$. Then for any prime $p$ with $p\mid |T|$, we have $p^d\mid d!$, and by Proposition~\ref{lmsylown}, $d<\frac{d}{p-1}$, a contradiction.

We have just shown that $N$ has one orbit, that is, $N$ is transitive on $V(\Ga)$.
If $N \cap G=1$, Lemma~\ref{lm-homo} implies that $G \cong K_v/N_v$ is solvable, a contradiction. Therefore, $G \leq N$, and by the arbitrariness of $N$, $X$ has only one minimal normal subgroup, that is, $\soc(X)=N$.

Since $G$ is not normal in $X$, we have $G<N$, and hence $N_v\not=1$ as $\Ga$ is $G$-vertex-transitive. Clearly, we may chose $v$ such that $N_v^{\Ga(v)}\not=1$. Since $\Ga$ has prime valency and $N_v^{\Ga(v)}\unlhd X_v^{\Ga(v)}$,  $N_v^{\Ga(v)}$ is transitive on $\Ga(v)$, that is, $\Ga$ is $N$-arc-transitive.

Recall that $N=T_1\times T_2\times\cdots\times T_d\cong T^d$. Suppose $d\geq 2$. If $T_1$ is transitive, then by Lemma~\ref{lm-homo}, $T_2\times\cdots\times T_d\cong N_v/(T_1)_v$ is solvable, a contradiction. Thus, $T_1$ has at least three orbits, and hence $|G|\geq 3|T_1|$. In particular, $G\cap T_1=1$. By the simplicity of $G$, the conjugate action of $G$ on $T_1$ is trivial or faithful. If it is trivial then $GT_1=G\times T_1$, and by Lemma~\ref{lm-homo}, $T_1\cong (GT_1)_v/G_v$ is solvable, a contradiction. Thus, the conjugate action of $G$ on $T_1$ is faithful and hence we may assume $G\leq \Aut(T_1)$. Since $|G|\geq 3|T_1|=3|\Inn(T_1)|$, we have $G\cap \Inn(T_1)=1$ and hence $G=G/(G\cap \Inn(T_1))\cong G\Inn(T_1)/\Inn(T_1)\leq \Aut(T_1)/\Inn(T_1)=\Out(T_1)$, which is impossible because $\Out(T_1)$ is solvable. Thus, $\soc(X)=N=T$ is a non-abelian simple group.

By the Frattini argument, $T=GT_v$. Then the triple $(T,G,T_v)$ can be read out from \cite{Li-Xia}, where $T_v$ is a group given in Proposition \ref{solvable-stabilizer}. Note that $p\geq 11$.

By \cite[Proposition 4.2]{Li-Xia}, $T$ cannot be any exceptional group of Lie type.

Assume that $T=\A_n$. By \cite[Proposition 4.3]{Li-Xia}, one of the following occurs:
\begin{itemize}
\item[(a)] $G=\A_{n-1},T=\A_n$ with $n \geq 6$ and $T_v$ is transitive on $\{1,2,\cdots, n\}$, or
\item[(b)] $G=\A_{n-2}, T=\A_n$ with $n=q^f$ for some prime $q$, and $T_v \leq \AGammaL(1,q^f)$ is 2-homogeneous on $\{1,2,\cdots, n\}$.
\end{itemize}

If (b) occurs, then $T_v$ is primitive on  $\{1,2,3,\cdots,q^f\}$ because it is 2-homogeneous. By Proposition~\ref{solvable-stabilizer}, $T_v$ has a normal subgroup $\ZZ_p$, and by the primitivity of $T_v$, $\ZZ_p$ is transitive and so regular on  $\{1,2,3,\cdots,q^f\}$. It follows $q^f=p$ and $T_v \leq \AGL(1,p)=\ZZ_p:\ZZ_{p-1}$. Moreover, since $|T_v|=\frac{|T||G_v|}{|G|}\geq \frac{|T|}{|G|}=p(p-1)$, we have that $T_v = \AGL(1,p)=\ZZ_p:\ZZ_{p-1}$. Thus, $\A_p$ contains a cyclic subgroup $\ZZ_{p-1}$, which is impossible because $\ZZ_{p-1}$ contains odd permutations on $\{1,2,3,\cdots,p\}$. It follows that $T=\A_n$, $G=\A_{n-1}$ and $T_v$ is transitive on the $n$ points, which is the case~(1) of the lemma.

Assume that $T$ is a sporadic simple group. By \cite[Proposition 4.4]{Li-Xia},$G=\M_{22}$, $T=\M_{23}$, and $T_v=\ZZ_{23}$ or $\ZZ_{23}:\ZZ_{11}$.  Suppose on the contrary that $T_v=\ZZ_{23}:\ZZ_{11}$. We may let $T_{uv}=\ZZ_{11}$ for $u \in \Ga(v)$. Since $\Ga$ is $T$-arc-transitive, there is an element $g \in T$ interchanging $u$ and $v$, and hence $T_{uv}^g=T_{u^gv^g}=T_{uv}$, that is, $g \in N_T(T_{uv})$. A computation with {\sc Magma}~\cite{Magma} shows that there is only one conjugate class of $\ZZ_{11}$ in $\M_{23}$, and the normalizer of $\ZZ_{11}$ in $\M_{23}$ is $\ZZ_{11}:\ZZ_5$. Thus, $g \in \ZZ_{11}:\ZZ_5$ has odd order, which is impossible because $g$ interchanges $u$ and $v$. It follows that $T_v=\ZZ_{23}$, which is the case~(2) of the lemma.

Assume that $T$ is a classical simple group of Lie type. Note that $T=GT_v$, $G$ is non-abelian simple and $T_v$ is solvable.
Let $H$ is a maximal subgroup subject to that $T_v\leq H$ and $H$ is solvable. Then $T=GH$, and  $(T,G,H)$ is listed in \cite[Table 1.1 and Table 1.2]{Li-Xia}. Clearly, $|T:G|\div |T_v| \div |H|$. For an integer $m$ and a prime $r$, we use $m_r$ to denote the largest $r$-power dividing $m$.

By Proposition~\ref{solvable-stabilizer}, $T_v=\ZZ_k \times (\ZZ_p: \ZZ_{\ell})$ with $k \div \ell \div p-1$, where $\ZZ_p: \ZZ_{\ell} \leq \AGL(1,p)$. Let $P$ and $Q$ be the maximal normal $r$-subgroup of $T_v$ and $H$ respectively. Then $Q \cap T_v \leq P$, and since $T_v/(T_v \cap Q) \cong QT_v/Q \leq H/Q$, we have $|T_v|_r \leq |T_v \cap Q|\cdot|H/Q|_r \leq |P||H/Q|_r$.
Clearly, $|T_v|_p=p$ and hence $|T:G|_p \leq  p$.

Suppose that $r \neq p$ and $r \div |T_v|$. If $P$ is not contained in $\ZZ_k$, then $1 \neq P\ZZ_k/\ZZ_k  \unlhd T_v/\ZZ_k\cong \ZZ_p:\ZZ_\ell $, which is impossible because $\ZZ_p$ is the unique minimal normal subgroup of $\ZZ_p:\ZZ_\ell $. Therefore $P \leq \ZZ_{k}$. It follows from $k \div \ell $ that $|P|^2 \leq |T_v|_r$, and from $|T_v|_r \leq |P||H/Q|_r$ that $|P|\leq |H/Q|_r$. Thus, $|T:G|_r \leq |T_v|_r \leq (|H/Q|_r)^2 $.

Since $G$ is a non-abelian simple group, we may exclude Row 1 of \cite[Table 1.1]{Li-Xia} and Rows 7-10, 17 and 21 of \cite[Table 1.2 ]{Li-Xia}, and since $p\geq 11$ and $p\div |H|$, we may exclude Rows 6, 11-13, 16-20, 22 and 24-27 of \cite[Table 1.2]{Li-Xia}. The remaining cases are Rows 2-9 of \cite[Table 1.1]{Li-Xia}, and Rows 1-5, 14, 15, 23 and 28 of \cite[Table 1.2]{Li-Xia}.

In what follows we write $q=r^f$ for some prime $r$ and positive integer $f$.

For Row 2 of \cite[Table 1.1]{Li-Xia}, $T=\PSL(4,q)$, $G=\PSp(4,q)$, and $H=q^3:\frac{q^3-1}{(4,q-1)}.3$. By \cite[Table 5.1A]{K-Lie}, $q^2 \div |T:G|$. Thus $r \neq p$. Note that $|H/Q|_r=1$ or $3$. Since $ r^{2f}=q^2 \leq |T_v|_r \leq |H/Q|_r^2$, we have $r=3$ and $f=1$, that is, $q=3$. This is impossible because a computation with {\sc Magma} shows $T=\PSL(4,3)$ has no factorization $T=GH$.

For Row 3 of \cite[Table 1.1]{Li-Xia}, $T=\PSp(2m,q)$, $G=\Omega^{-}(2m,q)$, and $H=q^{m(m+1)/2}:(q^{m}-1).m$ with $m \geq 2$ and $q$ even. Then $r=2$.
By \cite[Table 5.1A]{K-Lie}, $q^m \div |T:G|$, implying $r\not=p$. Furthermore, $r^{fm}=q^m \leq |T_v|_r \leq |H/Q|_r^2= m_r^2$. It follows $r^{m_2} \leq r^{fm} \leq m_2^2$, and this holds if and only if $m_2=2$ or $4$. If $m_2=2$, then $2^{fm} \leq  m_2^2=4$ implies $f=1,m=2$, which is impossible because $\PSp(4,2) \cong \S_6$ is not a simple group. If $m_2=4$, then $2^{fm} \leq m_2^2=16$ implies that $m=4$ and $f=1$. In this case, $|H|=2^{12} \cdot 15$, contradicting that $p \div |H|$ with $p\geq 11$.

For Rows 4 and 5 of \cite[Table 1.1]{Li-Xia}, $T=\PSp(4,q)$, $G=\PSp(2,q^2)$, and $H=q^{3}:\frac{q^{2}-1}{(2,q-1)}.2$. By \cite[Table 5.1A]{K-Lie}, $q^2 \div |T:G|$, and so $r \neq p$. Note that $|H/Q|_r=1$ or $2$. Since $ r^{2f}=q^2 \leq |T_v|_r \leq |H/Q|_r^2$, we have that $r=2$ and $f=1$. This is impossible because $T=\PSp(4,2) \cong \S_6$ is not simple.

For Row 6 of \cite[Table 1.1]{Li-Xia}, $T=\PSU(2m,q)$, $G=\SU(2m-1,q)$, and $H=q^{m^2}:\frac{q^{2m}-1}{q+1(2m,q+1)}.m$ with $m \geq 2$. By \cite[Table 5.1A]{K-Lie}, $q^{2m-1}=r^{(2m-1)f} \div |T:G|$ and $r \neq p$. Thus $r^{(2m-1)f}=q^{2m-1} \leq |T_v|_r \leq H/Q|_r^2= m_r^2$, implying $r^{2m_r-1}\leq m_r^2$, which is impossible.

For Row 7 of \cite[Table 1.1]{Li-Xia}, $T=\POmega(2m+1,q)$, $G=\Omega^{-}(2m,q)$, and $H=(q^{m(m-1)/2}.q^m):\frac{q^{m}-1}{2}.m$ with $m \geq 3$ and $q$ odd. Then $r,m_r \geq 3$. By \cite[Table 5.1A]{K-Lie}, $q^m=2^{fm} \div |T:G|$ and hence $r \neq p$. Then $r^{fm}=q^m \leq |T_v|_r \leq |H/Q|_r^2= m_r^2$, and so $r^{m_r}\leq m_r^2$, which is impossible.

For Row 8 of \cite[Table 1.1]{Li-Xia}, $T=\POmega^{+}(2m,q)$, $G =\Omega(2m-1,q)$, and $H=q^{m(m-1)/2}:\frac{q^{m}-1}{(4,q^m-1)}.m$ with $m \geq 5$. By \cite[Table 5.1A]{K-Lie}, $q^{m-1}=r^{f(m-1)} \div |T:G|$ and $r\neq p$. Then $r^{f(m-1)}=q^{m-1} \leq |T_v|_r \leq |H/Q|_r^2= m_r^2$. Note that the inequality $2^{x}> x^2$ always holds for $x \geq 5$. Thus $m_r \leq 4$. Since $r^{f(m-1)} \leq m_r^2$ and $m \geq 5$, we have that $r=2,m_r=4$ and $m=5$, which is impossible because $m_r=5_2=1$.

For Row 9 of \cite[Table 1.1]{Li-Xia}, $T=\POmega^{+}(8,q)$, $G =\Omega(7,q)$, and $H=q^{6}:\frac{q^{4}-1}{(4,q^4-1)}.4$. By \cite[Table 5.1A]{K-Lie}, $q^{3}=r^{3f} \div |T:G|$, and $r \neq p$. Then $r^{3f}=q^3 \leq |T_v|_r \leq |H/Q|_r^2= (4_r)^2$, implying $r=2$ and $f=1$. In this case, $|H|=2^8\cdot 15$, contradicting $p \div |H|$ with $p\geq 11$.

For Row 14 of \cite[Table 1.2 ]{Li-Xia}, $T=\PSp(4,11)$, $H=11_{+}^{1+2}:10.\A_4$, and $G=\PSL(2,11^2)$. By \cite[Table 5.1A]{K-Lie}, $11^{2} \div |T:G|\div |T_v|$ and hence $p \neq 11$, which is impossible because $p$ is the largest prime divisor of $|T_v|$. Similarly, we may exclude Row 15 of \cite[Table 1.2 ]{Li-Xia}, because $T=\PSp(4,23)$, $H=23_{+}^{1+2}:22.\S_4$, $G=\PSL(2,23^2)$, and  $23^{2} \div |T:G| \div |T_v|$ by \cite[Table 5.1A]{K-Lie}.

For Row 23 of \cite[Table 1.2 ]{Li-Xia}, $T=\Omega(7,3)$, $H=3^{3+3}:13:3$ and $G=\Sp(6,2)$. Then $p=13$, and since $|T_v|=pk\ell$ with $k\div \ell\div (p-1)$, we have $3^5\nmid |T_v|$. However, $|T:G|=|\Omega(7,3)|/|\Sp(6,2)|=13 \cdot 3^5$ implies $3^5\div |T_v|$, a contradiction. Similarly, we may exclude Row 28 of \cite[Table 1.2 ]{Li-Xia} because $T=\POmega^{+}(8,3)$, $H=3^6:(3^3:13:3)$ or $3^{3+6}:13:3$, $G=\Omega^{+}(8,2)$ and $|T:G|=13 \cdot 3^7$.

For Rows 1-5 of \cite[Table 1.2 ]{Li-Xia}, by {\sc Magma} we obtain the following:

\begin{itemize}

\item[(a)] $( G,T,T_v)=(\A_5,\PSL(2,11),\ZZ_{11})$,
\item[(b)] $( G,T,T_v)=(\A_5,\PSL(2,11),\ZZ_{11}:\ZZ_{5})$,
\item[(c)] $( G,T,T_v)=(\A_5,\PSL(2,19),\ZZ_{19}:\ZZ_{9})$,
\item[(d)] $( G,T,T_v)=(\A_5,\PSL(2,29),\ZZ_{29}:\ZZ_{7})$,
\item[(e)] $( G,T,T_v)=(\A_5,\PSL(2,29),\ZZ_{29}:\ZZ_{14})$,
\item[(f)] $( G,T,T_v)=(\A_5,\PSL(2,59),\ZZ_{59}:\ZZ_{29})$.
\end{itemize}
For case (b), $|V(\Ga)|=|T:T_v|=12$ and hence $\Ga$ is a complete graph of order $12$, contradicting that $\Aut(\Ga)_v$ is solvable. Similarly, cases (c),(e) and (f) cannot occur because $\Ga$ is a complete graph of order $20$, $30$ or $60$, respectively. Thus, we have (a) or (d), which is the case~(3) or (4) of the lemma.

For cases (2)-(4), it is easy to see that $G_v=G\cap T_v=1$. Since $G$ is transitive, it is  regular, that is, $\Ga$ is Cayley graph of $G$. \qed

\begin{lemma} \label{lmBlhd1}  Let $G$, $\Ga$, $p$ and $v$ be as given in Assumption and further assume that $G$ is regular on $V(\Ga)$. Then $\rad(\Aut(\Ga))$ has at least three orbits on $V(\Ga)$, and if $\rad(\Aut(\Ga))G\unlhd \Aut(\Ga)$  then $\rad(\Aut(\Ga))G=\rad(\Aut(\Ga)) \times G$.
\end{lemma}

\proof Set $A=\Aut(\Ga)$, $R=\rad(A)$ and $B=RG$. If $R$ is transitive on $V(\Ga)$, then Lemma \ref{lm-homo} implies that $G \cong B_v/R_v$ is solvable, a contradiction. Since $G$ is transitive, $\Ga$ is not bipartite and hence $R$ has at least three orbits. Assume that $B\unlhd A$. To finish the proof, it suffices to show $B=R\times G$. This is clearly true for $R=1$.

Assume $R\not=1$. Then $R\cap G=1$. Since $G$ is regular, $B_v\not=1$, and since $B\unlhd A$, $\Ga$ is $B$-arc-transitive. By Proposition~\ref{solvable-stabilizer}, $B_v$ has a normal sylow $p$-subgroup $\ZZ_p$, and $|B_v|=pm$ with $(p,m)=1$. Note that $RG=B=GB_v$. Again by the regularity of $G$, we have $|B_v|=|R|=pm$. Let $R_p$ be a Sylow $p$-subgroup of $R$. We claim $R_p\unlhd B$.

Suppose to the contrary that $R_p\ntrianglelefteq B$. Since $R\unlhd B$ is solvable, by the Jordan-Holder Theorem, $B$ has a normal series: $1\unlhd R_1\unlhd R_2\unlhd  \cdots \unlhd R\unlhd B$ such that $R_1\unlhd B$, $R_2\unlhd B$, $R_2/R_1 \cong \ZZ_p$ and $R_1\not=1$. Since $(p,m)=1$, we have $p\notdiv |R_1|$. Note that $R_2/R_1\unlhd B/R_1$ and $GR_1/R_1\cong G/(G\cap R_1)=G$. Since $R_2/R_1\cong\ZZ_p$, the conjugate action of $GR_1/R_1$ on $R_2/R_1$ must be trivial by the simplicity of $G$. It follows that $GR_2/R_1= R_2/R_1 \times GR_1/R_1$, and hence, $GR_1/R_1\unlhd GR_2/R_1$, forcing $GR_1\unlhd GR_2$. Since $p\mid |R_2|$, $GR_2$ is arc-transitive on $\Ga$, and hence $GR_1$ is also arc-transitive because $|(GR_1)_v|=|R_1|\not=1$. It follows $p\mid |R_1|$, a contradiction. Thus, $R_p\unlhd B$, as claimed.

Let $C=C_B(R_p)$. Since $R_p\unlhd B$ and $R_p\cong\ZZ_p$, the conjugate action of $G$ on $R_p$ is trivial and so $R_pG=R_p\times G$. It follows that $G\leq C$ and $C=C \cap B=C \cap (RG)=(C \cap R)G$.
Clearly, $R_p\leq C\cap R$ and hence $R_p$ is a Sylow $p$-subgroup of $C\cap R$. This implies that $C\cap R=R_p\times L$ where $L$ is a $p'$-subgroup of $C \cap R$, and in particular, $L$ is characteristic in $C\cap R$ and so normal in $B$.  Thus, $C=(R_p \times L)G=R_p \times LG$, and therefore, $LG\unlhd C$. Note that $C$ is arc-transitive because  $G\leq C$ and $R_p\leq C$. If $L\not=1$ then $(LG)_v\not=1$ and then $LG\unlhd C$ implies that $LG$ is arc-transitive. This means that $p\div |(LG)_v|$, and since $LG=G(LG)_v$, we have $|(LG)_v|=|L|$ and $p\div |L|$, which is impossible. It follows that $L=1$ and $C=R_p \times G$. Furthermore, $G\unlhd B$ and so $B=R \times G$. \qed

\medskip

\noindent {\bf Proof of Theorem~\ref{th1}:} Let $G$, $\Ga$, $p$ and $v$ as given in Assumption and further let $G$ be regular on $V(\Ga)$. Write $A=\Aut(\Ga)$, $R=\rad(A)$ and $B=RG$.
Then $R\cap G=1$ and $B/R\cong G$. By the Frattini argument, $B=GR=GB_v$, and so $|R|=|B_v|$.

Assume $R=1$. By Lemma~\ref{lm-R=1-trans}, either $\soc(A)=G$, or $\Ga$ is $\soc(A)$-arc-transitive and $G<\soc(A)$ with $(G,\soc(A))=(\A_{n-1},\A_n),(\M_{22},\M_{23})$, $(\A_5,\PSL(2,11))$ or $(\A_5,\PSL(2,29))$.

Assume $R\neq 1$. By Lemma~\ref{lmBlhd1}, $R$ has at least three orbits, and by Proposition~\ref{lm-arc}, the quotient graph $\Ga_R$ has valency $p$ with $A/R$-arc-transitive and $B/R$-vertex-transitive. Moreover, $(A/R)_\Delta \cong A_v$ is solvable for any $\Delta \in V(\Ga_R)$. Write $I/R=\soc(A/R)$. Since $B/R\cong G$, Lemma~\ref{lm-R=1-trans} implies that either $B/R=I/R \unlhd A/R$, or $\Ga_R$ is $I/R$-arc-transitive with $B/R<I/R$ and $(B/R,I/R)=(\A_{n-1},\A_n)$ with $(I/R)_\Delta$ being transitive on $\{1,2,\cdots,n\}$, or $(B/R,I/R,(I/R)_\Delta)=(\M_{22},\M_{23},\ZZ_{23}),(\A_{5},\PSL(2,11),\ZZ_{11}) $ or $( \A_{5},\PSL(2,29),\ZZ_{29}:\ZZ_{7})$.

\medskip
\noindent {\bf Case 1}: $B/R=I/R \unlhd A/R$.

In this case, $B=GR\unlhd A$, and by Lemma~\ref{lmBlhd1}, $B=G\times R$. It follows that $G$ is characteristic in $B$, and hence $G\unlhd A$.

\medskip
\noindent {\bf Case 2}: $\Ga_R$ is $I/R$-arc-transitive with $B/R<I/R$ and $(B/R,I/R)=(\A_{n-1},\A_n)$ with $(I/R)_\Delta$ being transitive on $\{1,2,\cdots,n\}$, or $(B/R,I/R,(I/R)_\Delta)=(\M_{22},\M_{23},\ZZ_{23})$, $(\A_{5},\PSL(2,11),\ZZ_{11}) $ or $( \A_{5},\PSL(2,29),\ZZ_{29}:\ZZ_{7})$.

Let $(B/R,I/R,(I/R)_\Delta)=(\M_{22},\M_{23},\ZZ_{23})$, $(\A_{5},\PSL(2,11),\ZZ_{11})$ or $( \A_{5}, \PSL(2,29),\ZZ_{29}:\ZZ_{7})$. By Lemma \ref{lm-R=1-trans}, $\Ga$ is a Cayley graph on $G B/R\cong G$. Since $\Ga$ is a Cayley graph on $G$, we have that $|V(\Ga)|=|V(\Ga_R)|$, which contradicts the assumption $R \neq 1$. Thus  $(B/R,I/R)=(\A_{n-1},\A_n)$ with $(I/R)_\Delta$ being transitive on $\{1,2,\cdots,n\}$.

First we claim $B=R \times G$. Suppose to the contrary that $B \neq R \times G$. Since $R$ is solvable, there exists a series of normal subgroups of $B$: $R_0 = 1<R_1<\dots<R_s=B$ such that $R_i \lhd B$ and $R_{i + 1}/ R_{ i}$ is an elementary abelian group for each $0 \leq i \leq s-1$. Since $RG \neq R \times G$, there exists $0\leq j \leq s-1$ such that $GR_i = G \times R_i$ for any $0 \leq i \leq j$, but $GR_{j+1} \neq G\times R_{j+1}$.

Write $R_{j+1}/R_j=\ZZ_r^f$ for some prime $r$ and positive integer $f$. Note that $G\cap R_i=1$ for $0\leq i\leq s$ and so $R_{i+1}G/R_i\cong G$ for $0 \leq i \leq s-1$. In particular, the conjugate action of $R_{j+1}G/R_j$ on $R_{j+1}/R_j$ is trivial or faithful. If it is trivial, then $R_{j+1}G/R_j=(R_{j+1}/R_j)(R_{j}G/R_j)=R_{j+1}/R_j \times R_{j}G/R_j$, implying $R_{j}G \lhd R_{j+1}G$, and since $GR_j = G \times R_j$, we have $G\unlhd R_{j+1}G$ and $GR_{j+1} = G \times R_{j+1}$, a contradiction. It follows that the conjugate action of $R_{j+1}G/R_j$ on $R_{j+1}/R_j$ is faithful, and we may assume $G\leq \GL(f,r)$.

Recall that $|B_v|=|R|$ and $R_{j+1}/R_j=\ZZ_r^f$. Then $r^f\div |B_v|$, and since $\Ga_R$ is $I/R$-arc-transitive, $\Ga$ is $I$-arc-transitive and  Proposition \ref{lm-arc} implies $I_v \cong (I/R)_\Delta$. Since $B/R<I/R$, we have $|B_v|\div |I_v|$ and so $r^f\div |(I/R)_\Delta|$. If $r=p$ then Proposition \ref{solvable-stabilizer} implies $r^2\nmid |(I/R)_\Delta|$ and hence $G\leq \GL(1,p)$, a contradiction. It follows $r\not=p$, and again by Proposition \ref{solvable-stabilizer}, $r^f\div (p-1)^2$.

Now $B/R=\A_{n-1}\leq \GL(f,r)$. By assumption, $p\geq 11$. Since $(I/R)_\Delta$ contains a normal subgroup $\mz_p$, we have $p\div n$ and so $n-1 \geq 11-1=10$. By Proposition \ref{lm-An}, $f \geq (n-1)-2 \geq p-3$ and so $(p-1)^2 \geq r^f \geq 2^{p-3}$. This is impossible because the function $f(x)=2^{x-3}-(x-1)^2>0$ always holds for $x \geq 11$. This completes the proof of the claim, and hence $B=R \times G$.

\vskip 0.2cm
Set $C=C_I(R)$. Then $G\leq C$, $C\unlhd I$ and $C\cap R\leq Z(C)$. Recall that $I/R=\A_n $ or $\M_{23}$. Since $G\cong (R\times G)/R\leq CR/R \unlhd I/R$, we have $I=CR$, and since $Z(C)/(C\cap R) \unlhd C/C\cap R\cong CR/R=I/R$, we have $C\cap R=Z(C)$ and $C/Z(C)\cong I/R$.   Furthermore, $C'/(C'\cap Z(C))\cong C'Z(C)/Z(C)=(C/Z(C))'=C/Z(C)\cong  I/R$, and so  $Z(C')=C'\cap Z(C)$, $C=C'Z(C)$ and $C'/Z(C')\cong I/R$. It follows $C'=(C'Z(C))'=C''$, and hence $C'$ is a covering group of $I/R$.

Suppose $Z(C')\not=1$. Then Proposition \ref{lmcoveran} implies that $Z(C')=\ZZ_2$ and $C'\cong 2.\A_n$. Since $G\leq C$ and $C/C'$ is abelian, we have $G\leq C'$. So $G\times Z(C') \cong \A_{n-1} \times \ZZ_2$ is a subgroup of $C' \cong 2.\A_n$, which is impossible by Proposition \ref{lmcoveran}.

Thus, $Z(C')=1$. It follows $C'\cong I/R$. Since $G<C$ and $C/C'$ is abelian, we have $G<C'\unlhd I$, and since $|I|=|I/R||R|=|C'||R|$ and $C'\cap R=1$, we have $I=C'\times R$. Since $C^{'}$ is a nonabelian simple group, $C'$ is characteristic in $I$, and hence $C'\unlhd A$ because $I\unlhd A$. Since $G$ is regular on $\Ga$ and $G<C'\unlhd I$, $C'$ has non-trivial stabilizer, and hence $\Ga$ is $C^{'}$-arc-transitive on $\Ga$. Note that $C'\cong I/R=\A_n$.

\vskip 0.2cm
Summing up, we have proved that either $G\unlhd A$, or $A$ has a normal arc-transitive subgroup $T$ such that $G<T$ and $(G,T)=(\A_5,\PSL(2,11)), (\A_5,\PSL(2,29)), (\M_{22},\M_{23})$ or $(\A_{n-1},\A_n)$ (for $R=1$, $T=\soc(A)$, and for $R\not=1$, $T=C'$). Let  $(G,T)=(\A_{n-1},\A_n)$. Since $G$ is regular, $|T_v|=n$, and by Proposition~\ref{solvable-stabilizer}, $n=pk\ell$ with $k\div \ell\div (p-1)$. To finish the proof, we are left to show that $k$ and $\ell$ have the same parity.

Suppose to the contrary that $k$ and $\ell$ has different parity. Then $k$ is odd and $\ell$ is even as $k\div \ell$. Since $(G,T)=(\A_{n-1},\A_n)$, we have $|T:G|=n$ and $T$ can be viewed as the alternating permutation group by the well-known right multiplication action of $T$ on the set $[T:G]$ of all right cosets of $G$ in $T$, still denoted by $\A_n$. By the regularity of $G$ on $\Ga$, $T=GT_v$ and $G\cap T_v=1$, which implies that $T_v\leq \A_n$ is a regular permutation group on $[T:G]$. By Proposition~\ref{solvable-stabilizer}, $T_v=\mz_k\times(\mz_p: \mz_\ell)$, and so $T_v$ has a cyclic group $\mz_\ell$ with odd index $|T_v:\mz_\ell|=pk$. Let $\mz_\ell=\langle a\rangle$. Since $T_v$ is regular, $a$ is a product of $pk$ $\ell$-cycles on $[T:G]$ in its distinct cycle decomposition, so an odd permutation as $\ell$ is even and $kp$ is odd, which is impossible because $T_v\leq \A_n$. This completes the proof.     \qed

\section{Proof of Theorem~\ref{th2}}\label{nonnormalexample}

The goal of this section is to prove Theorem \ref{th2}. To do that, we first describe a widely known construction for vertex-transitive and symmetric graphs, part of which
is attributed to Sabidussi~\cite{Sabi}.

Let $G$ be a group, $H$ a subgroup of $G$, and $D$ a union of some double cosets of $H$ in $G$ such that $H\not\subseteq D$ and $D^{-1}=D$. Then the \emph{coset graph} $\Ga=\Cos(G,H,D)$ is defined as the graph with vertex-set $[G\!:\!H]$, the set of all right cosets of $H$ in $G$, and edge-set $E(\Ga)=\{\{Hg,Hxg\} : g\in G, x\in D\}$. This graph is regular with valency $|D|/|H|$, and is connected if and only if $G=\langle D,H\rangle$, that is, if and only if $G$ is generated by $D$ and $H$. The group $G$ acts vertex-transitively on $\Ga$ by right multiplication. More precisely, for $g\in G$, the permutation $\hat{g}_{H}: Hx\mapsto Hxg$, $x\in G$, on $[G:H]$ is an automorphism of $\Aut(\Ga)$, and $\hat{G}_{H}:=\{\hat{g}_{H}\ |\ g\in G\}$ is a transitive subgroup of $\Aut(\Ga)$.
The map $g\mapsto \hat{g}_H$, $g\in G$, is a homomorphism from $G$ to $S_{[G:H]}$, the well-known coset action of $G$ on $H$, and the kernel of this coset action is $H_G=\bigcap_{g\in G}H^g$, the largest normal subgroup of $G$ contained in $H$. It follows that $G/H_G\cong \hat{G}_{H}$. Furthermore, $\Ga$ is $\hat{G}_H$-arc-transitive if and only if $D$ consists of just one double coset $HaH$. If $H_G=1$, we say that $H$ is {\em core-free} in $G$, and in this case, $G\cong \hat{G}_{H}$.

If $H=1$, denote $\Cos(G,H,D)$ and $\hat{G}_H$ by $\Cay(G,D)$ and $\hat{G}$, respectively. In this case, $\hat{G}$ is the right regular representation of $G$, and it is regular on the vertex set of $\Cay(G,D)$. By definition, $\Cay(G,D)$ is Cayley graph of $\hat{G}$, and for short,  $\Cay(G,D)$ is also called a Cayley graph of $G$ with respect to $D$.

Conversely, suppose $\Ga$ is any graph on which the group $G$ acts faithfully and vertex-transitively. Then it is easy to show that $\Ga$ is isomorphic to the coset graph $\Cos(G,H,D)$, where $H=G_{v}$ is the stabiliser in $G$ of the vertex $v\in V(\Ga)$, and $D$ is a union of double cosets of $H$, consisting of all elements of $G$ taking $v$ to one of its neighbours. Then $H\not\subseteq D$ and $D^{-1}=D$.
Moreover, if $G$ is arc-transitive on $\Ga$ and $g $ is an element of $G$ that swaps $v$ with one of its neighbours,  then $g^2 \in H$ and $D=HgH$, and the valency of $\Ga$ is $|D|/|H|=|H:H\cap H^g|$. Also $a$ can be chosen as a $2$-element in $G$. In particular, if $L\leq G$ is regular on vertex set of $\Ga$, then $\Ga$ is also isomorphic to $\Cay(L,S)$, where $S$ consists of all elements of $L$ taking $v$ to one of its neighbours with $S^{-1}=S$, and by the regularity, we have $S=D\cap L$. Thus, we have the following proposition.

\begin{proposition}\label{doublecosetgraph}
Let $\Ga$ be a $G$-vertex-transitive graph and $L$ be a regular subgroup of $G$. Then $\Ga\cong \Cos(G,H,D)\cong \Cay(L,S)$ with $S=L\cap D$, where $H=G_{v}$ for $v\in V(\Ga)$, $D$ is a union of double cosets of $H$, consisting of all elements of $G$ taking $v$ to one of its neighbours, and $S$ consists of all elements of $L$ taking $v$ to one of its neighbours. Moreover, $\Ga$ be $G$-arc-transitive if and only if $G$ has a $2$-element $g$ such that $D=HgH$, and in this case, $\Ga$ has valency $|H:H\cap H^g|$.
\end{proposition}

Let $\Ga=\Cos(G,H,D)$ be a coset graph. We set $\Aut(G, H, D)=\{\alpha \in \Aut(G)\ |\ H^{\alpha}=H, D^{\alpha}=D\}$. For any $\alpha\in \Aut(G, H, D)$, the permutation $\alpha_{H}: Hx\mapsto Hx^{\alpha}$, $x\in G$, on $[G:H]$ is an automorphism of $\Ga$, and the map $\a\mapsto \a_H$ is a natural action of $\Aut(G, H, D)$ on $V(\Ga)$. It follows that $\Aut(G, H, D)/K\cong \Aut(G, H, D)_{H}$, where $\Aut(G, H, D)_{H}=\{\alpha_{H}|\ \alpha\in \Aut(G,H,D)\}$ and $K$ is the kernel of the action. Furthermore,  $\Aut(G,H,D)_H\leq \Aut(\Ga) $.
For $h\in H$, let $\tilde{h}$ be the inner automorphism of $G$ induced by $h$, that is,  $\tilde{h}:g\mapsto h^{-1}gh$, $g\in G$. Then $\tilde{H}:=\{\tilde{h} \ | \ h\in H \}\leq \Aut(G,H,D)$ and hence $\tilde{H}_{H}:=\{\tilde{h}_{H}\ |\ h \in H\}$ is a subgroup of $\Aut(G, H, D)_{H}$.

The following proposition was proved by Wang, Feng and Zhou~\cite[Lemma~2.10]{WFZ}, which is important for computing automorphism groups of coset graphs.

\begin{proposition}\label{CosetgraphAuto}
Let $G$ be a finite group, $H$ a core-free subgroup of
$G$ and $D$ a union of several double-cosets $HgH$
such that $H\nsubseteq D$ and $D=D^{-1}$. Let $\Ga=\Cos(G,H,D)$ and $A=\Aut(\Ga)$. Then $\hat{G}_H\cong G$, $\Aut(G,H,D)_{H}\cong \Aut(G,H,D)$, $\tilde{H}_{H}\cong \tilde{H}$, and  $\N_{A}(\hat{G}_{H})=\hat{G}_{H} \Aut(G, H, D)_{H}$ with $\hat{G}_{H}\cap
\Aut(G, H, D)_{H}=\tilde{H}_{H}$.
\end{proposition}

Now we are ready to prove Theorem~\ref{th2} and this follows from  Lemmas~\ref{th2-lemma1}-\ref{th2-lemma4}.

\medskip
Let $x,y,t$ be permutations in $\S_{11}$ as following:

\vskip -0.4cm
\begin{eqnarray*}
  x&=&(1, 11, 8, 3, 6, 9, 4, 10, 2, 7, 5)\\
  y&=&(2, 10, 6)(3, 11, 4)(7, 8, 9) \\
  t& =& (2, 5)(3, 9)(6, 11)(8, 10)
\end{eqnarray*}

Let $T=\langle x,t\rangle$, $H=\langle x \rangle$, $G=\langle y,t\rangle$. Define
$$\Ga=\Cos(T,H,HtH).$$

Then a computation with {\sc Magma} \cite{Magma} shows that $T \cong \PSL(2,11)$, $H \cong \ZZ_{11}$, $|H \cap H^t|=1$, and $G \cong \A_5$. By Proposition~\ref{doublecosetgraph}, $\Ga$ has valency $11$ and $T$ acts arc-transitively on $\Ga$. Since $11 \notdiv |G|$, $G$ acts semiregularly on $V(\Ga)$, and since $|G|=|V(\Ga)|$, $G$ is regular on $V(\Ga)$. It follows that $\Ga$ is a non-normal Cayley group of $\A_5$ with $\PSL(2,11)$-arc-transitive. A direct computation with {\sc Magma}  shows that $\Aut(\Ga) \cong \PGL(2,11)$ and this implies the following lemma.

\begin{lemma}\label{th2-lemma1}
 There exists an $11$-valent symmetric Cayley graph $\Ga$ of $A_5$ such that $\Aut(\Ga) \cong \PGL(2,11)$. In particular, $\Aut(\Ga)_v$ is solvable for $v\in V(\Ga)$.
\end{lemma}

\medskip
Let $x,y,t,z$ be permutations in $\S_{30}$ as following:

\vskip -0.6cm
\begin{eqnarray*}
  x&=&(1, 21, 10, 9, 22, 28, 13, 15, 30, 6, 19, 18, 7, 27, 23, 4, 25, 17, 20, 2, 12, 29, 16, 26, 8,11, \\
   & & 3, 24, 5) \\
  y&=& (1, 24, 9)(2, 6, 5)(3, 27, 21)(4, 12, 20)(7, 25, 26)(8, 10, 13)(11, 14, 16)(15, 30, 23)(17, \\
  & &  28,29)(18, 22, 19) \\
  t& =& (1, 3)(2, 10)(4, 11)(5, 19)(6, 24)(7, 16)(8, 17)(9, 28)(12, 27)(13, 20)(14, 22)(15, 26)\\
   & &(18, 30)(21, 23)\\
 z& =& (2, 18, 23, 10, 29, 9, 17)(3, 7, 19, 20, 4, 24, 30)(5, 22, 27, 13, 28, 6, 16)(8, 12, 15, 21, 11,\\
  & &   25, 26)
\end{eqnarray*}

Let $T=\langle x,t\rangle$, $H=\langle x,z \rangle$, $G=\langle y,t\rangle$. Define
$$\Ga=\Cos(T,H,HtH).$$

Then a computation with {\sc Magma} \cite{Magma} shows that $T \cong \PSL(2,29)$, $H \cong \ZZ_{29}:\ZZ_7$, $|H \cap H^t|=7$, and $G \cong \A_5$. By Proposition~\ref{doublecosetgraph}, $\Ga$ has valency $29$ and $T$ acts arc-transitively on $\Ga$. Since $29 \notdiv |G|$, $G$ acts semiregularly on $V(\Ga)$, and since $|G|=|V(\Ga)|$, $G$ is regular on $V(\Ga)$. It follows that $\Ga$ is a non-normal Cayley group of $\A_5$ with $\PSL(2,29)$-arc-transitive. A direct computation with {\sc Magma}  shows that $\Aut(\Ga) \cong \PGL(2,29)$ and this implies the following lemma.

\begin{lemma}\label{th2-lemma2}
There exists a $29$-valent symmetric Cayley graph $\Ga$ of $A_5$ such that $\Aut(\Ga) \cong \PGL(2,29)$. In particular, $\Aut(\Ga)_v$ is solvable for $v\in V(\Ga)$.
\end{lemma}

\medskip
Let $x,y,t$ be permutations in $\S_{23}$ as following:

\vskip -0.6cm
\begin{eqnarray*}
  x&=&(1, 4, 6, 7, 2, 19, 3, 11, 9, 20, 13, 23, 16, 8, 21, 5, 14, 22, 18, 15, 17, 10, 12)\\
  y&=&(1, 14, 6, 5, 9, 2, 10, 3, 15, 13, 11)(4, 22, 16, 19, 17, 8, 21, 7, 12, 18, 23) \\
  t& =& (1, 17)(3, 9)(5, 18)(6, 13)(7, 12)(10, 19)(14, 22)(21, 23)
\end{eqnarray*}

Let $T=\langle x,t\rangle$, $H=\langle x \rangle$, $G=\langle y,t\rangle$. Define
$$\Ga=\Cos(T,H,HtH).$$

\begin{lemma}\label{th2-lemma3} The above graph $\Ga$ is $23$-valent symmetric Cayley graph of $M_{22}$ and $\Aut(\Ga)=(\hat{\M}_{23})_H\cong \M_{23}$. In particular, $\Aut(\Ga)_v$ is solvable for $v\in V(\Ga)$.
\end{lemma}

\proof A computation with {\sc Magma} \cite{Magma} shows that $T \cong \M_{23}$, $H \cong \ZZ_{23}$, $|H \cap H^t|=1$, and $G \cong \M_{22}$. By Proposition~\ref{doublecosetgraph}, $\Ga$ has valency $23$ and $T$ acts arc-transitively on $\Ga$. Since $23 \notdiv |G|$, $G$ acts semiregularly on $V(\Ga)$, and since $|G|=|V(\Ga)|$, $G$ is regular on $V(\Ga)$. It follows that $\Ga$ is a non-normal Cayley group of $\M_{22}$ with $\M_{23}$-arc-transitive. However, we cannot compute $\Aut(\Ga)$ with {\sc Magma} because $|V(\Ga)|$ is too large. By Proposition~\ref{doublecosetgraph}, we may let $\Ga=\Cay(G,S)$ with $S=G\cap HtH$. Write $A=\Aut(\Ga)$. By {\sc Magma}, $S=\{s_i\ |\ 1\leq i\leq 23 \}$, where
\vskip 0.1cm
$s_1=(1, 14, 6, 5, 9, 2, 10, 3, 15, 13, 11)(4, 22, 16, 19, 17, 8, 21, 7, 12, 18, 23)$,

$s_2=(1, 11, 13, 15, 3, 10, 2, 9, 5, 6, 14)(4, 23, 18, 12, 7, 21, 8, 17, 19, 16, 22)$,

 $s_3=(1, 15, 5, 2, 12, 18, 16, 14, 21, 13, 7)(3, 6, 4, 22, 8, 19, 10, 17, 9, 23, 11)$,

$s_4=(1, 7, 13, 21, 14, 16, 18, 12, 2, 5, 15)(3, 11, 23, 9, 17, 10, 19, 8, 22, 4, 6)$,

 $s_5=(1, 9, 14)(2, 19, 5, 4, 22, 12)(3, 21, 6)(7, 23, 15, 11, 8, 18)(10, 13)(16, 17)$,

$s_6=(1, 14, 9)(2, 12, 22, 4, 5, 19)(3, 6, 21)(7, 18, 8, 11, 15, 23)(10, 13)(16, 17)$,

$s_7=(1, 4, 3)(2, 6)(5, 8, 7, 10, 14, 21)(9, 12, 17, 22, 16, 13)(11, 19, 23)(15, 18)$,

$s_8=(1, 3, 4)(2, 6)(5, 21, 14, 10, 7, 8)(9, 13, 16, 22, 17, 12)(11, 23, 19)(15, 18)$,

$s_9=(1, 12)(2, 19, 3)(4, 6, 18, 5, 8, 10)(7, 11, 23, 16, 14, 22)(9, 13)(15, 17, 21)$,

$s_{10}=   (1, 12)(2, 3, 19)(4, 10, 8, 5, 18, 6)(7, 22, 14, 16, 23, 11)(9, 13)(15, 21, 17)$,

$s_{11}=(1, 7, 3, 16, 12)(2, 11, 23, 22, 14)(4, 15, 5, 18, 10)(6, 9, 13, 8, 17)$,

 $s_{12}=(1, 12, 16, 3, 7)(2, 14, 22, 23, 11)(4, 10, 18, 5, 15)(6, 17, 8, 13, 9)$,

  $s_{13}=  (3, 16, 23, 12, 6)(4, 11, 22, 18, 10)(5, 17, 7, 19, 9)(8, 14, 15, 21, 13)$,

 $s_{14}=(3, 6, 12, 23, 16)(4, 10, 18, 22, 11)(5, 9, 19, 7, 17)(8, 13, 21, 15, 14)$,

  $s_{15}=(1, 15, 12, 6, 19)(2, 11, 13, 14, 7)(3, 16, 21, 22, 4)(5, 10, 17, 9, 23)$,

   $s_{16}=(1, 19, 6, 12, 15)(2, 7, 14, 13, 11)(3, 4, 22, 21, 16)(5, 23, 9, 17, 10)$

   $s_{17}=(1, 7)(3, 8)(4, 6)(9, 19)(11, 23)(12, 15)(13, 18)(14, 21)$,

    $s_{18}=(2, 6)(3, 10)(4, 22)(8, 16)(11, 13)(12, 18)(14, 15)(21, 23)$,

    $s_{19}=(1, 11)(2, 16)(4, 19)(6, 12)(8, 14)(9, 13)(15, 18)(17, 22)$,

    $s_{20}=(1, 17)(3, 9)(5, 18)(6, 13)(7, 12)(10, 19)(14, 22)(21, 23)$,

    $s_{21}=(1, 15)(5, 16)(6, 18)(7, 19)(8, 21)(9, 23)(11, 12)(17, 22)$,

     $s_{22}=(1, 17)(2, 9)(5, 11)(6, 19)(7, 13)(8, 23)(10, 12)(14, 15)$,

    $s_{23}= (1, 5)(2, 4)(3, 11)(8, 13)(9, 19)(10, 15)(14, 16)(18, 23)$.
\vskip 0.1cm
Let $1$ be the identity in $G$. Then $1\in V(\Ga)$. Suppose to the contrary that $A_1$ is nonsolvable. By Proposition~\ref{localgloblestabilizers}, the restriction $A_1^{\Ga(1)}$ of $A_1$ on the neighbourhood $\Ga(1)$ of $1$ in $\Ga$ is nonsolvable, and since $\Ga$ has prime valency, the Burnside Theorem (also see \cite[Theorem  3.5B]{Dixon}) implies that $A_1^{\Ga(1)}$ is $2$-transitive on $\Ga(1)$. This turns that there exists a $5$-cycle passing though $1$ and any two vertices in $S$ because $(1,s_{11},s_{11}^2,s_{11}^3,s_{11}^4)$
is a $5$-cycle in $\Ga$. In particular, there is a $5$-cycle passing through $1$, $s_1$ and $s_2=s_1^{-1}$, and hence $s_1^2\in S^3=\{s_{i_1}s_{i_2}s_{i_2}\ |\ s_{i_1},s_{i_2},s_{i_2}\in S\}$, but this is not true by {\sc Magma}~\cite{Magma}. Thus, $A_1$ is solvable.

Now we let $\Ga=\Cos(T,H,HtH)$ and $D=HtH$. Since $A$ has solvable stabilizer, Theorem~\ref{th1} implies that $\hat{T}=\hat{\M}_{23} \unlhd A$. Note that $H$ is core-free in $T$. By Proposition~\ref{CosetgraphAuto},
$A=\hat{T}_{H} \Aut(T, H, D)_{H}$ with $\hat{T}_{H}\cap
\Aut(T, H, D)_{H}=\tilde{H}_{H}$, where $\hat{T}_H\cong T$, $\Aut(T,H,D)_{H}\cong \Aut(T,H,D)$ and $\tilde{H}_{H}\cong \tilde{H}$. To prove $A=\hat{T}_{H}$, it suffices to show that $\Aut(T,H,D)=\tilde{H}$.

Suppose to the contrary that $\a\in \Aut(T,H,D)$, but $\a\not\in \tilde{H}$. By \cite[Table 5.1.C]{K-Lie}, $\Out(\M_{23})=1$, that is, $\Aut(\M_{23})=\Inn(\M_{23})$. Thus, $\a$ is an automorphism of $T$ induced by an element of $b\in T$ by conjugation, namely $g^\a=g^b$ for $g\in T$. Since $\a\in \Aut(T,H,D)$, we have $H^b=H$ and $D^b=D$, and since  $\a\not\in \tilde{H}$, we have $b\not\in H$. It follows that $H\langle b\rangle$ is a subgroup of $T$ containing $H$, and by {\tt Atlas}~\cite{CCNPW}, $H\langle b\rangle\cong \mz_{23}:\mz_{11}$. Since $\tilde{H}\leq \Aut(T,H,D)$, we may choose $b$ such that $b$ has order $11$, and by {\sc Magma}, we may let  $b=(2, 14, 18, 7, 16, 6, 9, 20, 8, 3, 4)(5, 21, 13, 22, 12, 15, 11, 19, 17, 23, 10)$ because $H=\langle x\rangle$ with
$x=(1, 4, 6, 7, 2, 19, 3, 11, 9, 20, 13, 23, 16, 8, 21, 5, 14, 22, 18, 15, 17, 10, 12)$. However, $D^b=(HtH)^b\not=HtH$ by {\sc Magma}, a contradiction. Thus, $A=\hat{T}_{H}\cong \M_{23}$. \qed

\medskip
 Let $p \geq 5$ be a prime, and let $x$, $t$ and $h$ be permutations in $\S_p$ as following:
 \vskip -0.4cm
\begin{eqnarray*}
x= ( 1,2,\cdots,p), & t = (1 ,2)( 3,4 ), & h=(2,p)(3,p-1)\cdots (\frac{p-1}{2},\frac{p+5}{2})(\frac{p+1}{2},\frac{p+3}{2}).
\end{eqnarray*}

Let $T = \langle x , t \rangle$ and $H= \langle x \rangle$. By \cite{FMW}, $T=\A_p$, $H\cong \ZZ_{p}$ and $|H \cap H^t|=1$. Define
$$\Ga^p=\Cos(\A_p,H,HtH).$$

\begin{lemma}\label{th2-lemma4} The above graph $\Ga^p$ is a $p$-valent symmetric Cayley graph of $A_{p-1}$ such that
$\Aut(\Ga^p)\cong\S_p$ for $p\equiv 3 (\mod 4)$ and $\Aut(\Ga^p)\cong\A_p \times \ZZ_2$ for $p\equiv 1 (\mod 4)$. In particular, $\Aut(\Ga)_v$ is solvable for $v\in V(\Ga)$.
\end{lemma}

\proof By Proposition~\ref{doublecosetgraph}, $\Ga^p$ has valency $p$ and $\A_p$ acts arc-transitively on $\Ga^p$, with vertex stabilizer isomorphic to $\mz_p$. Let $\A_{p-1}$ be the subgroup of $\A_p$ fixing the point $p$. Since $p \notdiv |\A_{p-1}|$, $\A_{p-1}$ acts semiregularly on $V(\Ga^p)$, and since $|\A_{p-1}|=|V(\Ga^p)|$, $\A_{p-1}$ is regular on $V(\Ga^p)$. It follows that $\Ga^p$ is a non-normal Cayley group of $\A_{p-1}$ with $\A_p$-arc-transitive.

By Proposition~\ref{doublecosetgraph}, we may let $\Ga^p=\Cay(\A_{p-1},S)$, where $S=\A_{p-1} \cap HtH$.
For $p=5$ or $p=7$, a computing with {\sc Magma} shows that $\Aut(\Ga^5)\cong \A_5\times \mz_2$ and $\Aut(\Ga^7)\cong \S_7$. Write $A=\Aut(\Ga)$. We may assume $p\geq 11$.

\vskip 0.2cm
\noindent {\bf Claim:} $A$ has solvable stabilizer.

 Recall that $x=( 1,2,\cdots,p)$, $t = (1 ,2)( 3,4 )$ and $H= \langle x \rangle$. Let $x^{-i}tx^j \in S=HtH\cap \A_{p-1}$ for $i,j \in\mz_p$. Then $p=p^{x^{-i}tx^j}=p^{x^{-i}tx^ix^{j-i}}$. Note that
$x^{-i}tx^i=(1^{x^i},2^{x^i})(3^{x^i},4^{x^i})$, and if $j-i\not=0$ then $x^{j-i}$ is a $p$-cycle. For $0\leq i\leq p-5$, $p=p^{x^{-i}tx^ix^{j-i}}=p^{x^{j-i}}$ implies $j=i$. Furthermore, For $i=p-4,p-3,p-2$ or $p-1$,  $p=p^{x^{-i}tx^ix^{j-i}}$ implies that $j=i+1,i-1,i+1$ or $i-1$, respectively. Thus, we may set $S=\{s_1,s_2,\cdots,s_p\}$, where

\vskip 0.1cm
\noindent $s_{i+1}=x^{-i}tx^i=(1+i,2+i)(3+i,4+i)$ for $0 \leq i \leq p-5$,

\noindent $s_{p-2}=x^{-(p-3)}tx^{p-4}=(1,p-1,p-3,\cdots,3,2), \ \ \ s_{p-3}=x^{-(p-4)}tx^{p-3}=(s_{p-2})^{-1}$,

\noindent $s_p=x^{-(p-1)}tx^{p-2}=(1,p-1,p-2, \cdots,4,3), \ \ \ s_{p-1}=x^{-(p-2)}tx^{p-1}=s_p^{-1}$.

\vskip 0.1cm
For $z\in \A_p$, denote by $o(z)$ the order of $z$ and by $supp(z)$  the support of $z$, that is, the number of points moving by $z$. Then $o(s_i)=2$ and $supp(s_i)=4$ for $1\leq i\leq p-4$, and $o(s_i)=supp(s_i)=p-2$ for $p-3\leq i\leq p$.

To prove the Claim, it suffices to show that $A_1$ is solvable. We argue by contradiction and we suppose that $A_1$ is nonsolvable. Note that $\Ga^p=\Cay(\A_{p-1},S)$ and $\Ga^p(1)=S$.

By Propostion~\ref{localgloblestabilizers}, $A_1^{\Ga^p(1)}$ is nonsolvable, and the Burnside Theorem implies that $A_1$ is $2$-transitive on $\Ga^p(1)$. Note that $p\geq 11$. Since $s_1=(1,2)(3,4)$ commutes with $s_5=(5,6)(7,8)$, there is a $4$-cycle passing through $1,s_1$ and $s_5$. By the $2$-transitivity of $A_1$ on $\Ga^p(1)$, there exists a $4$-cycle through $1,s_p$ and $s_{p-1}=s_p^{-1}$, and this implies $|Ss_p\cap Ss_p^{-1}|\geq 2$. Thus, $|Ss_p^{-2}\cap S|\geq 2$.

Let $S_1=\{s_i\ |\ 1\leq i\leq p-4\}$ and $S_2=\{s_{p-2},s_{p-2}^{-1},s_p,s_p^{-1}\}$. Then $S=S_1\cup S_2$ and $S_1\cap S_2=\emptyset$. Since $s_p^{-1}$ is a $(p-2)$-cycle in $\A_p$ and $p-2$ is odd, $s_p^{-2}$ is also a $(p-2)$-cycle, implying $supp(s_p^{-2})=p-2$. Since $supp(s_i)=4$ for each $1\leq i\leq p-4$, we have $supp(s_is_p^{-2})\geq p-6\geq 5$, and  $s_is_p^{-2}$ cannot be any involution in $S$. Thus, $|S_1s_p^{-2}\cap S_1|=0$.

Note that $S_2s_p^{-2}=\{s_p^{-1},s_p^{-3},s_{p-2}s_p^{-2},s_{p-2}^{-1}s_p^{-2}\}$. Then  $|S_2s_p^{-2}\cap S_2|=1$ by a simple checking one by one. If $|S_2s_p^{-2}\cap S_1|\not=0$, then $z^2=1$ for some $z\in S_2s_p^{-2}$, and we have $s_p^{-2}=1$ or $s_p^{-6}=1$ or $(s_{p-2}s_p^{-2})^2=1$ or $(s_{p-2}^{-1}s_p^{-2})^2=1$, of which all are impossible because all these elements cannot fix $1$. Thus, $|S_2s_p^{-2}\cap S_1|=0$. Similarly, $|S_2s_p^2\cap S_1|=0$.

Recall that $|Ss_p^{-2}\cap S|\geq 2$. Since $|S_2s_p^{-2}\cap S_2|=1$ and $|S_2s_p^{-2}\cap S_1|=0$, we have $|S_1s_p^{-2}\cap S|=1$, and since $|S_1s_p^{-2}\cap S_1|=0$, we have $|S_1s_p^{-2}\cap S_2|=1$. It follows $|S_2s_p^2\cap S_1|=1$, a contradiction. Thus, $A_1$ is solvable, as claimed.

\vskip 0.2cm
From now on, we write $\Ga^p=\Cos(T,H,HtH)$. Clearly, $H$ is core-free in $T$. By Claim, $A=\Aut(\Ga^p)$ has solvable stabilizer. By Theorem~\ref{th1}, $\hat{T}_H$ is normal in $A$, and by Proposition \ref{CosetgraphAuto}, $A=N_A(\hat{T}_H)=\hat{T}_H\Aut(T,H,HtH)_H$ with $\hat{T}_H \cap \Aut(T,H,HtH)_H=\tilde{H}_H$. Furthermore, $\hat{T}_H\cong T$, $\tilde{H}_H\cong H$ and $\Aut(T,H,HtH)_H\cong \Aut(T,H,HtH)=\{ \alpha \in \Aut(T)|H^\alpha=H,(HtH)^\alpha =HtH\}$.

Let $x^itx^j\in HtH$ for some $i,j\in \mz_p$. If $i+j=0$, then $x^itx^j=(1+j,2+j)(3+j,4+j)$ and $supp(x^itx^j)=4$. If $i+j \neq 0$, then $x^{i+j}$ is a $p$-cycle and $supp(x^itx^j)=supp(x^{i+j}x^{-j}tx^{j})\geq p-4>4$ because $supp(x^{-j}tx^{j})=4$. Thus, $I:=\{x^{-i}tx^i\ | i\in\mz_p\}$ consists of all elements in $HtH$ whose supports are 4.

Now we consider $\Aut(T,H,HtH)$. Let $\b\in \Aut(T,H,HtH)$. Then $\b\in \Aut(T)=\Aut(\A_p)\cong \S_p$, and $\b$ is an automorphism of $T$ induced by some $b\in \S_p$ by conjugation, that is, $t^\b=t^b$ for any $t\in T$. Since $(HtH)^\b=(HtH)^b=HtH$, we have $I^\b=I$, and in particular, $supp(yz)=supp(y^\b z^\b)$ for any $y,z\in I$. It is easy to see that for any $x^{-i}tx^i, x^{-j}tx^j\in I$, $supp(x^{-i}tx^ix^{-j}tx^j)=5$ if and only if $j=i+1$ or $i-1$. In fact, if $j=i+2$ or $i-2$ then $supp(x^{-i}tx^i x^{-j}tx^j)=4$, if $j=i+3$ or $i-3$ then $supp(x^{-i}tx^i x^{-j}tx^j)=7$, and if $|i-j|\geq 4$ then  $supp(x^{-i}tx^i x^{-j}tx^j)=8$.

Let $\Sigma$ be a graph with $I$ as vertex set and with $y,z\in I$ adjacent if and only if $supp(yz)=5$. By the above paragraph, $\Sigma$ is a cycle of length $p$, and $\b$ induces an automorphism of $\Sigma$. Thus, $\Aut(T,H,HtH)$ acts on $I$, and since $\Sigma$ is a $p$-cycle, $\Aut(T,H,HtH)/K\leq D_{2p}$, where $K$ is the kernel of this action. Let $\g\in K$, and suppose $\g$ is induced by $c\in \S_p$ by conjugation. Then $\g$ fixes each element in $I$, that is, $(x^{-i}tx^i)^c=x^{-i}tx^i$ for each $i\in \mz_p$. Since $x^{-i}tx^i=(1^{x^i},2^{x^i})(3^{x^i},4^{x^i})$ and $x^{-(i+3)}tx^{i+3}=(4^{x^i},5^{x^i})(6^{x^i},7^{x^i})$, $c$ fixes $\{1^{x^i},2^{x^i}, 3^{x^i},4^{x^i}\}$ and $\{4^{x^i},5^{x^i}, 6^{x^i},7^{x^i}\}$ setwise, and hence fixes $4^{x^i}=\{1^{x^i},2^{x^i}, 3^{x^i},4^{x^i})\}\cap \{4^{x^i},5^{x^i}, 6^{x^i},7^{x^i})\}$ for each $i\in\mz_p$. It follows that $c$ fixes $\{1,2,\cdots,n\}$ pointwise, implying $K=1$. Thus, $|\Aut(T,H,HtH)|\leq |\Aut(\Sigma)|=2p$.

Recall that $h=(2,p)(3,p-1)\cdots (\frac{p-1}{2},\frac{p+5}{2})(\frac{p+1}{2},\frac{p+3}{2})$. For $p=1\ \mod 4$, $h$ is an even permutation and $h\in \A_p$, and for $p=3\ \mod 4$, $h$ is an odd permutation and $h\in \S_p$, but $h\not\in\A_p$. Since $x=( 1,2,\cdots,p)$, we have $x^h=x^{-1}$ and so $H^h=H$, and since $t^h=(1^h,2^h)(3^h,4^h)=(1,p)(p-1,p-2)=x^{-(p-3)}tx^{p-3}\in I\subset HtH$, we have $(HtH)^h=HtH$. Clearly, $H^x=H$ and $(HtH)^x=H$. For any $z\in\S_p$, denote by $\tilde{z}$ the induced automorphism of $\A_p$ by $z$ by conjugation. Then $\tilde{x},\tilde{h}\in \Aut(T,H,HtH)$ and $\langle \tilde{x},\tilde{h}\rangle\cong D_{2p}$. Since $|\Aut(T,H,HtH)|\leq 2p$, we have $\Aut(T,H,HtH)=\langle \tilde{x},\tilde{h}\rangle\cong D_{2p}$.

Recall that $\tilde{x}_H: Hg\mapsto Hg^x$ for $g\in \A_p$, and $\tilde{h}_H: Hg\mapsto Hg^h$ for $g\in \A_p$, are automorphisms of $\Ga^p$, and $\tilde{H}_H=\langle\tilde{x}_H \rangle$. Since $\Aut(T,H,HtH)\cong \Aut(T,H,HtH)_H$, we have $\Aut(T,H,HtH)_H=\langle \tilde{x}_H,\tilde{h}_H\rangle=\tilde{H}_H:\tilde{h}_H\cong \D_{2p}$, and since $\hat{T}_H\cap \Aut(T,H,HtH)_H=\tilde{H}_H$ and $A=\hat{T}_H\Aut(T,H,HtH)_H$, we have $|A:\hat{T}_H|=2$ and hence $A=\hat{T}_H:  \langle\tilde{h}_H\rangle$.

Set $C=C_A(\hat{T}_H)$, the centralizer of $\hat{T}_H$ in $A$. Since $\hat{T}_H\cong \A_p$, we have $C\cap \hat{T}_H=1$, and since $A=\hat{T}_H: \langle\tilde{h}_H\rangle$, we have $C=1$ or $C\cong \mz_2$. For the former, $A\cong \S_p$ by the N/C Theorem, and for the latter, $A=\hat{T}_H\times C\cong \A_p\times\mz_2$. To finish the proof, we only need to prove that $C\cong \mz_2$ if and only if  $p=1\ \mod 4$.

Assume $C\cong \mz_2$. Since $A=\hat{T}_H\rtimes \langle\tilde{h}_H\rangle$, we can let $C=\langle  \hat{y}_H\tilde{h}_H \rangle $ for some $y \in T$. This implies that for any $z,g\in T$, we have $(Hz)^{\hat{y}_H\tilde{h}_H\hat{g}_H}=(Hz)^{\hat{g}_H\hat{y}_H\tilde{h}_H}$, that is, $H(zy)^hg=H(zgy)^h$, implying $Hhzyhg=Hhzgyh $. Set $\ell=yhg(gyh)^{-1}$. Then $Hhz\ell(hz)^{-1}=H$, that is, $\ell\in H^{hz}=H^z$ for any $z\in \A_p$. This implies that $\ell\in \bigcap_{z\in \A_p} H^z$, and since $\bigcap_{z\in \A_p} H^z$ is the largest normal subgroup of of $\A_p$ contained in $H$, we have $\bigcap_{z\in \A_p} H^z=1$ and hence $\ell=1$. This means that $yhg=gyh$, and by the arbitrary of $g\in \A_p$, we have $yh\in C_{\A_p}(\S_p)=1$. It follows that $h=y \in \A_p$ and hence $p=1\ \mod 4$. On the other hand, if $p=1\ \mod 4$ then it is easy to check that $\hat{h}\tilde{h}\in C$. Thus, $C\cong \mz_2$ if and only if $p=1\ \mod 4$, as required. \qed

\end{document}